\documentclass[letter,twocolumn]{article}

\usepackage{graphicx}
\usepackage{times}
\usepackage{color}
\usepackage{fix2col}
\usepackage{amsmath,amsfonts,amssymb}

\newcommand{\QED}{$\Box$}
\newcommand{\eq}{\triangleq}

\DeclareMathOperator*{\argmin}{arg\,min}
\DeclareMathOperator*{\supp}{supp}

\DeclareMathOperator{\dez}{{\mathcal D}}
\DeclareMathOperator{\shrink}{{\mathcal S}}
\DeclareMathOperator{\sat}{sat}

\newcommand{\field}[1]{\mathbb{#1}}
\newcommand{\R}{\field{R}}

\newcommand{\U}{{\mathcal{U}}}

\newcommand{\diag}{\mathrm{diag}}
\newcommand{\blockdiag}{\mathrm{blockdiag}}
\newcommand{\vc}[1]{{\boldsymbol{#1}}}

\newtheorem{thm}{Theorem}

\newtheorem{lem}[thm]{Lemma}
\newtheorem{prop}[thm]{Proposition}
\newtheorem{defn}[thm]{Definition}
\newtheorem{rem}[thm]{Remark}
\newtheorem{problem}[thm]{Problem}

\title{\LARGE \bf
Maximum-Hands-Off Control and $L^1$ Optimality%
\thanks{
This research is supported in part by the JSPS Grant-in-Aid for Scientific Research (C) No.~24560543,
and also
Australian Research Council's
Discovery Projects funding scheme (project number DP0988601). 
}
}
\author{Masaaki Nagahara%
\thanks{
M. Nagahara is with
 Graduate School of Informatics, Kyoto
 University, Kyoto, 606-8501, 
 Japan; 
 email: {\tt nagahara@ieee.org}.
},
Daniel E. Quevedo%
\thanks{
D. E. Quevedo is with School of Electrical Engineering \& 
 Computer Science, The University of Newcastle, NSW
 2308, Australia;
 email: {\tt dquevedo@ieee.org}.
},
Dragan Ne\v{s}i\'{c}%
\thanks{
D. Ne\v{s}i\'{c} is with
 Department of Electrical and Electronic Engineering, 
 The University of Melbourne, Victoria 3010
 Australia; email: {\tt dnesic@unimelb.edu.au}
}
}
\date{}

\begin{document}
\maketitle
\thispagestyle{empty}
\pagestyle{empty}

\begin{abstract}
In this article, 
we propose a new paradigm of control,
called a maximum-hands-off control.
A hands-off control is defined as a control that has
a much shorter support than the horizon length.
The maximum-hands-off control is the minimum-support (or sparsest)
control among all admissible controls.
We first prove that
a solution to an $L^1$-optimal control problem gives a maximum-hands-off control, and vice versa.
This result rationalizes the use of $L^1$ optimality in computing a maximum-hands-off control.
The solution has in general the "bang-off-bang" property,
and hence the control may be discontinuous.
We then propose an $L^1$/$L^2$-optimal control to obtain a continuous hands-off control.
Examples are shown to illustrate the effectiveness of the proposed control method.
\end{abstract}

\section{Introduction}
In practical control systems, 
we often need to
minimize the control effort
so as to achieve control objectives
under limitations in equipment such as
actuators, sensors, and networks.
For example, 
the energy (or $L^2$-norm) of a control signal is minimized
to prevent engine overheating
or to reduce transmission cost
with a standard LQ (linear quadratic) control problem;
see e.g., \cite{AndMoo}.
Another example is the \emph{minimum-fuel} control,
discussed in e.g., \cite{Ath63,AthFal},
in which the total expenditure of fuel is minimized with
the $L^1$ norm of the control.

Alternatively, in some situations, the control effort can be dramatically reduced by
holding the control value \emph{exactly zero} over a time interval.
We call such control a \emph{hands-off control}.
A motivation for hands-off control is a stop-start system
in automobiles.
It is a hands-off control; it automatically shuts down 
the engine to avoid it idling for long periods of time.
By this, we can reduce CO or CO2 emissions as well as fuel consumption
\cite{Dun74}.
This strategy is also used in hybrid vehicles \cite{Cha07};
the internal combustion engine is stopped when
the vehicle is at a stop or the speed is lower than a preset threshold,
and the electric motor is alternatively used.
Thus hands-off control is also available for solving environmental problems.
Hands-off control is also
desirable for networked and embedded systems
since the communication channel is not used
during a period of zero-valued control.
This property is advantageous in particular for wireless communications
\cite{JeoJeo06,KonWonTsa09}.
In other words, hands-off control is the least \emph{attention}
in such periods.
From this point of view, hands-off control that maximizes the total time of no attention is somewhat related to the concept of minimum attention control
\cite{Bro97}.
Motivated by these applications, we propose 
a new paradigm of control, called \emph{maximum-hands-off control}
that maximizes the time interval over which the control is exactly zero.

The hands-off property is related to \emph{sparsity},
or the \emph{$L^0$ ``norm''} 
(the quotation marks indicate that this is not a norm;
see Section~\ref{sec:preliminaries} below)
of a signal,
defined by the total length of the intervals over which the signal
takes non-zero values. The maximum-hands-off control,
in other words, seeks 
the \emph{sparsest} (or $L^0$-optimal) control among all admissible controls.
This problem is however hard to solve since the cost function
is non-convex and discontinuous.
To overcome the difficulty, one can adopt
$L^1$ optimality as a convex approximation of the problem,
as often used in \emph{compressed sensing},
which has recently attracted significant attention
in signal processing;
see \cite{Ela,EldKut,HayNagTan13} for details.
Compressed sensing has shown by theory and experiments that 
sparse high-dimensional signals
can be reconstructed from incomplete measurements
by using $L^1$ optimization;
see e.g., \cite{Don06,Can06}.

Interestingly,  an $L^1$-optimal (or minimum-fuel) control has been known to have
such a sparsity property, traditionally called \emph{"bang-off-bang"}
as in \cite{AthFal}.
Although advantage has implicitly taken of the sparsity property for minimizing the fuel consumption,
there has been no research on the theoretical connection between sparsity and
$L^1$ optimality of the control.
In this article, we prove that 
a solution to an $L^1$-optimal control problem gives a maximum-hands-off control, and vice versa.
As a result,
the sparsest solution (i.e., the maximum-hands-off control)
can be obtained by solving an $L^1$-optimal control problem.

A consequence of our result is that
maximum hands-off control necessarily has a
"bang-off-bang" property;
the control abruptly changes
its values between $0$ and $\pm u_{\max}$
at switching times ($u_{\max}$ is the admissible maximum absolute value
of control).
In some applications, this feature should be avoided.
To make the control continuous in time,
we propose a new type of control, namely,
\emph{$L^1$/$L^2$-optimal control}.
We show that the $L^1$/$L^2$-optimal control is an intermediate control
between the maximum-hands-off (or $L^1$-optimal)
and the minimum energy (or $L^2$-optimal) controls,
in the sense that the $L^1$ and $L^2$ controls are
the limiting instances
of the $L^1$/$L^2$-optimal control.


The remainder of this article is organized as follows.
In Section~\ref{sec:preliminaries},
we give mathematical preliminaries for our subsequent discussion.
In Section~\ref{sec:problems},  we define two control problems:
maximum-hands-off control and $L^1$-optimal control.
In Section~\ref{sec:L1}, we briefly review $L^1$-optimal control.
Section~\ref{sec:main} gives the main theorem,
establishing the theoretical connection between the $L^1$-optimal control
and the maximum-hands-off one.
In Section~\ref{sec:L1L2}, we propose a mixed $L^1$/$L^2$-optimal control 
that gives a continuous hands-off control.
Section~\ref{sec:examples} presents control design examples
to illustrate the effectiveness of our method.
In Section~\ref{sec:conclusion}, we offer concluding remarks.

\section{Mathematical Preliminaries}
\label{sec:preliminaries}
For a vector $\vc{v}=[v_1,v_2,\dots,v_m]^\top$,
we define the $\ell^1$, $\ell^2$, and $\ell^\infty$ norms respectively by
\[
 \begin{split}
  &\|\vc{v}\|_1 \eq \sum_{i=1}^m |v_i|,\quad
  \|\vc{v}\|_2 \eq \sqrt{\sum_{i=1}^m |v_i|^2},\\
  &\|\vc{v}\|_\infty \eq \max_{i=1,\dots,m} |v_i|.
 \end{split}
\]

For a continuous-time signal $u(t)$
over a time interval $[0,T]$,
we define its $L^p$ norm ($p>0$)
by
\[
 \|u\|_{L^p} \eq \left(\int_0^T |u(t)|^p dt\right)^{1/p}.
\]
Note that if $p\in(0,1)$, then $\|\cdot\|_{L^p}$ is not a
norm (It fails to satisfy the triangle inequality.).
We define the support set of $u$, denoted by $\supp(u)$, by
the closure of the set
\[
\{t\in[0,T]: u(t)\neq0\}.
\]
Then we define
the $L^0$ ``norm'' of $u$ as
the length of its support, that is,
\[
 \|u\|_{L^0} \eq \mu\bigl(\supp(u)\bigr),
\] 
where $\mu$ is the Lebesgue measure on $\R$.
Note that the $L^0$ ``norm'' is not a norm since
it fails to satisfy the positive homogeneity, that is,
for any non-zero scalar $\alpha$ such that
$|\alpha|\neq 1$, we have
\[
 \|\alpha u\|_{L^0} = \|u\|_{L^0} \neq |\alpha| \|u\|_{L^0},
 \quad \forall u\neq 0.
\]
The notation $\|\cdot\|_{L^0}$ may be however
justified from the fact that
if $u$ is integrable on $[0,T]$, then
$u\in L^p$ for any $p\in(0,1)$ and
\[
 \lim_{p\rightarrow 0}\|u\|_{L^p}^p = \|u\|_{L^0},
\]
which is proved by using H\"{o}lder's inequality; see \cite{Tao} for details.
A continuous-time signal $u$ over $[0,T]$ is called
\emph{sparse}%
\footnote{This is analogous to the definition of the $\ell^0$ ``norm''
of a vector, which is defined by the number of non-zero elements.
When a vector has small $\ell^0$ ``norm,'' 
then it is also called sparse.
See \cite{Ela,EldKut,HayNagTan13} for details.}
if $\|u\|_{L^0}$ is ``much smaller'' than $T$.

For a function $\vc{f}=[f_1,\dots,f_n]^\top: \R^n\rightarrow\R^n$,
the Jacobian $\vc{f}'$ is defined by
\[
 \vc{f}'(\vc{x}) \eq 
  \begin{bmatrix}
   \frac{\partial f_1}{\partial x_1}&\dots&\frac{\partial f_1}{\partial x_n}\\
   \vdots & \ddots & \vdots\\
   \frac{\partial f_n}{\partial x_1}&\dots&\frac{\partial f_n}{\partial x_n}\\
 \end{bmatrix},
\]  
where $\vc{x}=[x_1,\dots,x_n]^\top$.

\section{Optimal Control Problems}
\label{sec:problems}
We here consider nonlinear plant models of the form
\begin{equation}
 \frac{d\vc{x}(t)}{dt} = \vc{f}\bigl(\vc{x}(t)\bigr) + \sum_{i=1}^m \vc{g}_i\bigl(\vc{x}(t)\bigr)u_i(t),
 \quad t\in[0,T],
 \label{eq:plant}
\end{equation}
where
$\vc{x}$ is the state,
$u_1,\dots,u_m$ are the control inputs,
$\vc{f}$ and $\vc{g}_i$
are functions on $\R^n$.
We assume that $\vc{f}(\vc{x})$, $\vc{g}_i(\vc{x})$,
and their Jacobians $\vc{f}'(\vc{x})$, $\vc{g}_i'(\vc{x})$
are continuous in $\vc{x}$.
We use the vector representation $\vc{u}\eq[u_1,\dots,u_m]^\top$.

The control $\{\vc{u}(t): t\in[0,T]\}$ is chosen to drive the state $\vc{x}(t)$
from a given initial state 
\begin{equation}
 \vc{x}(0)=\vc{x}_0,
 \label{eq:initial_state}
\end{equation} 
to the origin by a fixed final time $T>0$, that is,
\begin{equation}
 \vc{x}(T)=\vc{0}.
 \label{eq:final_state}
\end{equation}
Also, the control $\vc{u}(t)$ is constrained in magnitude by
\begin{equation}
 \|\vc{u}(t)\|_\infty \leq 1,\quad \forall t \in [0,T].
 \label{eq:input_constraint}
\end{equation}
We call a control $\{\vc{u}(t): t\in[0,T]\}$ \emph{admissible}
if it satisfies \eqref{eq:input_constraint}
and the resultant state $\vc{x}(t)$ from \eqref{eq:plant} satisfies boundary conditions
\eqref{eq:initial_state} and \eqref{eq:final_state}.
We denote by $\U$ the set of all admissible controls.

The \emph{maximum-hands-off control}
is a control that
maximizes the time interval over which the control $\vc{u}(t)$ is exactly zero.
In other words, we try to find the \emph{sparsest} control
among all admissible controls in $\U$.

We state the associated optimal control problem as follows:
\begin{problem}[Maximum-Hands-Off Control]
\label{prob:MHO}
Find an admissible control $\{\vc{u}(t): t\in[0,T]\}\in\U$ that minimizes
\begin{equation}
 J_0(\vc{u}) \eq \sum_{i=1}^m\lambda_i \|u_i\|_{L^0},
 \label{eq:J_MHO}
\end{equation}
where $\lambda_1>0,\dots,\lambda_m>0$ are given weights.
\end{problem}

On the other hand, if we replace $\|u_i\|_{L^0}$ in \eqref{eq:J_MHO}
with the $L^1$ norm $\|u_i\|_{L^1}$,
we obtain the following \emph{$L^1$-optimal control} problem,
also known as \emph{minimum fuel control}
discussed in e.g. \cite{Ath63,AthFal}.
\begin{problem}[$L^1$-Optimal Control]
\label{prob:L1}
Find an admissible control $\{\vc{u}(t): t\in[0,T]\}\in\U$ that minimizes
\begin{equation}
 J_1(\vc{u}) \eq \sum_{i=1}^m \lambda_i \|u_i\|_{L^1} = \int_0^T \sum_{i=1}^m \lambda_i |u_i(t)| dt,
 \label{eq:J_L1}
\end{equation}
where $\lambda_1>0,\dots,\lambda_m>0$ are given weights.
\end{problem}

\begin{rem}[Minimum time]
For the existence of the solution of both problems above,
the final time $T$ must be sufficiently large.
More precisely, $T$ must be larger than
the minimum time $T^\ast$ required to force
the initial state $\vc{x}_0$ to the origin.
$T^\ast$ is obtained by solving the minimum-time problem;
see \cite[Chap.~6]{AthFal} for details.
\end{rem}

\section{Review of $L^1$-Optimal Control}
\label{sec:L1}
Here we briefly review the $L^1$-optimal (or minimum-fuel) control problem (Problem~\ref{prob:L1})
based on the discussion in \cite[Sec.~6-13]{AthFal}.

\begin{figure}[tb]
\centering
\includegraphics[width=0.9\linewidth]{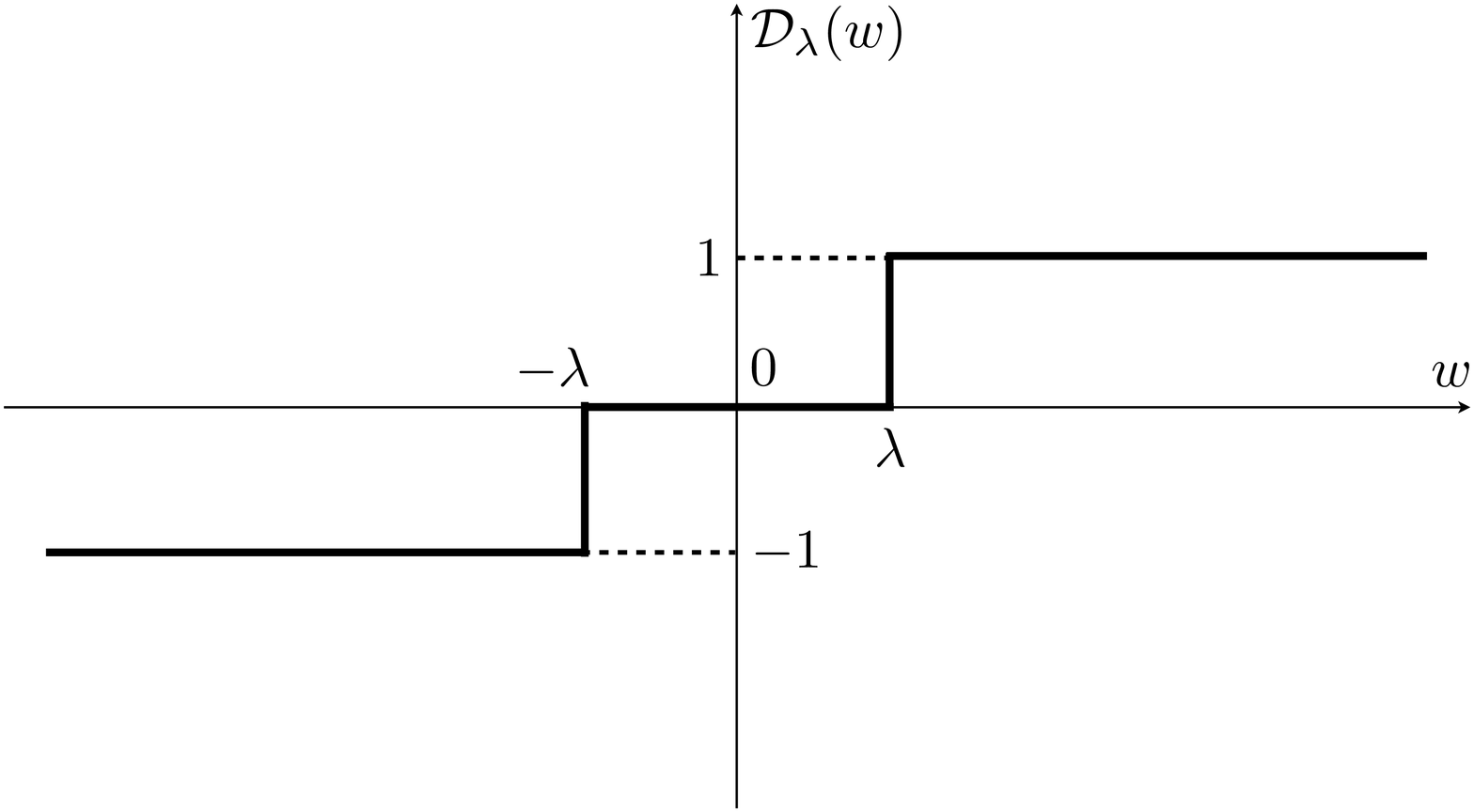}
\caption{Dead-zone function $\dez_\lambda(w)$}
\label{fig:dez}
\end{figure}

Let us first form the Hamiltonian function for the $L^1$-optimal control problem as
\begin{equation}
 H(\vc{x},\vc{p},\vc{u}) = \sum_{i=1}^m \lambda_i |u_i| 
  + \vc{p}^\top \biggl(\vc{f}\bigl(\vc{x}\bigr)
    +\sum_{i=1}^m \vc{g}_i(\vc{x})u_i\biggr),
 \label{eq:Hamiltonian_L1}
\end{equation}
where 
$\vc{p}$ is the costate (or adjoint) vector.
Assume that $\vc{u}^\ast=[u_1^\ast,\dots,u_m^\ast]^\top$
is an $L^1$-optimal control
and $\vc{x}^\ast$ is the resultant trajectory.
According to the minimum principle, there exists a costate $\vc{p}^\ast$ such that
the optimal control $\vc{u}^\ast$ satisfies
\[
 H(\vc{x}^\ast,\vc{p}^\ast,\vc{u}^\ast)\leq  H(\vc{x}^\ast,\vc{p}^\ast,\vc{u}),
\]
for all admissible $\vc{u}$.
The optimal state $\vc{x}^\ast$ and costate $\vc{p}^\ast$ 
satisfies the canonical equations
\[
 \begin{split}
  \frac{d\vc{x}^\ast(t)}{dt} &=  \vc{f}\bigl(\vc{x}^\ast(t)\bigr) + \sum_{i=1}^m
    \vc{g}_i\bigl(\vc{x}^\ast(t)\bigr)u_i^\ast(t),\\
  \frac{d\vc{p}^\ast(t)}{dt} &= - \vc{f}'\bigl(\vc{x}^\ast(t)\bigr)^\top \vc{p}^\ast(t)\\
  	&\qquad-\sum_{i=1}^m u_i^\ast(t) \vc{g}_i'\bigl(\vc{x}^\ast(t)\bigr)^\top \vc{p}^\ast(t),
 \end{split} 
\]
with boundary conditions
\[
 \vc{x}^\ast(0) = \vc{x}_0,\quad \vc{x}^\ast(T) = \vc{0}.
\]
The minimizer
$\vc{u}^\ast=[u_1^\ast,\dots,u_m^\ast]^\top$ of the Hamiltonian given in 
\eqref{eq:Hamiltonian_L1} is given by
\[
 u_i^\ast(t) = -\dez_{\lambda_i}\left(\vc{g}_i\bigl(\vc{x}^\ast(t)\bigr)^\top \vc{p}^\ast(t)\right),\quad t\in[0,T],
\] 
where $\dez_{\lambda}(\cdot):\R^n\rightarrow[-1,1]$ is the dead-zone function defined by
\begin{equation}
 \begin{split}
 \dez_{\lambda}(w)  &=
   \begin{cases} 
     -1,& \text{~if~} w<-\lambda,\\
     0, & \text{~if~} -\lambda<w<\lambda,\\
     1, & \text{~if~} \lambda<w,\\
   \end{cases}\\
  \dez_{\lambda}(w) &\in [-1,0], \text{~if~} w=-\lambda,\\
  \dez_{\lambda}(w) &\in [0,1], \text{~if~} w=\lambda.
 \end{split}
 \label{eq:dez}
\end{equation}
See Fig.~\ref{fig:dez} for the graph of $\dez_\lambda(\cdot)$.

If $\vc{g}_i(\vc{x}^\ast)^\top \vc{p}^\ast$ is equal to $-\lambda_i$ or $\lambda_i$ over
a non-zero time interval, say $[t_1,t_2]\subset[0,T]$,
$t_1<t_2$,
then the control $u_i$
(and hence $\vc{u}$)
over $[t_1,t_2]$
cannot be uniquely determined
by the minimum principle.
In this case, the interval $[t_1,t_2]$ is called a \emph{singular interval},
and a control problem that has at least one singular interval is called
\emph{singular}.
If there is no singular interval, the problem is called \emph{normal}:
\begin{defn}[Normality]
The $L^1$-optimal control problem stated in Problem~\ref{prob:L1} is said to be \emph{normal}
if the set 
\[
 {\mathcal T}_i \eq \{t\in [0,T]: |\lambda_i^{-1}\vc{g}_i(\vc{x}^\ast(t))^\top \vc{p}^\ast(t)|=1\}
\]
is countable for $i=1,\dots,m$.
If the problem is normal, the elements $t_1,t_2,\dots\in {\mathcal T}_i$
are called the \emph{switching times} for the control $u_i(t)$.
\hfill\QED
\end{defn}

If the problem is normal, the components
of the $L^1$-optimal control $\vc{u}^\ast(t)$
are piecewise constant and ternary,
taking values $\pm 1$ or $0$
at almost all $t\in[0,T]$.
This property, named "bang-off-bang,"
is key to connect the $L^1$-optimal control and
the maximum-hands-off control
as discussed in the next section.

\section{Maximum-Hands-Off Control and $L^1$-Optimal Control}
\label{sec:main}
In this section, we consider a theoretical relation
between
maximum-hands-off control (Problem \ref{prob:MHO})
and $L^1$-optimal control (Problem \ref{prob:L1}).
The theorem below rationalizes the $L^1$ optimality
in computing the maximum-hands-off control.
\begin{thm}
\label{thm:L1optimal}
Assume that the $L^1$-optimal control problem stated in Problem~\ref{prob:L1}
is normal and has at least one solution.
Let $\U_0^\ast$ and $\U_1^\ast$ be the sets of the optimal solutions
of Problem \ref{prob:MHO} ($L^0$-optimal control problem)
and Problem \ref{prob:L1} ($L^1$-optimal control problem)
respectively.
Then we have $\U_0^\ast = \U_1^\ast$.
\end{thm}

{\bf Proof.}
Let $\U$ be the set of all admissible controls for the
$L^1$-optimal control problem (Problem~\ref{prob:L1}).
By assumption, $\U_1^\ast$ is non-empty, and so is $\U$.
The set $\U$ is also the admissible control set
for the $L^0$-optimal control problem (Problem \ref{prob:MHO}),
and hence $\U_0^\ast \subset \U$.
We first show that $\U_0^\ast$ is non-empty, and then prove $\U_0^\ast = \U_1^\ast$.

First, for any $\vc{u}\in{\mathcal U}$, we have
\begin{equation}
 \begin{split}
  J_1(\vc{u}) &= \sum_{i=1}^m \lambda_i \int_0^T |u_i(t)| ~dt\\
  &= \sum_{i=1}^m \lambda_i \int_{\supp(u_i)} |u_i(t)| ~dt\\  
  &\leq \sum_{i=1}^m \lambda_i \int_{\supp(u_i)} 1 ~dt
  = J_0(\vc{u}).
 \end{split} 
\label{eq:proof1}
\end{equation}

Now take an arbitrary $\vc{u}^\ast_1 \in \U_1^\ast$.
Since the problem is normal by assumption,
each control $u_{1i}^\ast(t)$ in $\vc{u}_1^\ast(t)$ takes values $-1$, $0$, or $1$,
at almost all $t\in[0,T]$.
This implies that
\begin{equation}
 \begin{split}
 J_1(\vc{u}^\ast_1) &= \sum_{i=1}^m \lambda_i \int_0^T |u_{1i}^\ast(t)| ~dt\\
  &= \sum_{i=1}^m \lambda_i \int_{\supp(u_{1i}^\ast)} 1 ~dt
  = J_0(\vc{u}^\ast_1).
 \end{split}
 \label{eq:proof2}
\end{equation}
From \eqref{eq:proof1} and \eqref{eq:proof2},
$\vc{u}^\ast_1$ is a minimizer of $J_0$,
that is, $\vc{u}_1^\ast\in\U_0^\ast$.
Thus, $\U_0^\ast$ is non-empty and $\U_1^\ast \subset \U_0^\ast$.

Conversely, let $\vc{u}^\ast_0\in\U_0^\ast\subset\U$.
Take independently $\vc{u}^\ast_1\in\U_1^\ast\subset\U$.
From \eqref{eq:proof2} and the optimality of $\vc{u}^\ast_1$, we have
\begin{equation}
 J_0(\vc{u}^\ast_1) = J_1(\vc{u}^\ast_1) \leq J_1(\vc{u}^\ast_0).
 \label{eq:proof3}
\end{equation} 
On the other hand, from \eqref{eq:proof1} and the optimality of $\vc{u}^\ast_0$,
we have
\begin{equation}
 J_1(\vc{u}^\ast_0)\leq J_0(\vc{u}^\ast_0) \leq J_0(\vc{u}^\ast_1).
 \label{eq:proof4}
\end{equation}
It follows from \eqref{eq:proof3} and \eqref{eq:proof4} that
$J_1(\vc{u}^\ast_1)=J_1(\vc{u}^\ast_0)$,
and hence $\vc{u}^\ast_0$ achieves the minimum value of $J_1$.
That is, $\vc{u}_0^\ast\in\U_1^\ast$ and $\U_0^\ast\subset\U_1^\ast$.
\hfill\QED

Theorem~\ref{thm:L1optimal} suggests that
$L^1$ optimization can be used for 
the maximum-hands-off (or the sparsest) solution.
This is analogous to the situation in compressed sensing,
where $L^1$ optimality is often used to obtain the sparsest vector;
see \cite{Ela,EldKut,HayNagTan13} for details.


\section{$L^1$/$L^2$-Optimal Control}
\label{sec:L1L2}
In the previous section, we have shown that
the maximum-hands-off control problem can be solved
via $L^1$-optimal control.
From the "bang-off-bang" property of the $L^1$-optimal control,
the control changes its value at switching times \emph{discontinuously}.
This is undesirable for some applications in which the actuators cannot move
abruptly.
In this case, one may want to make the control \emph{continuous}.
For this purpose, we add a regularization term to the $L^1$ cost $J_1(\vc{u})$
defined in \eqref{eq:J_L1}.
More precisely, we consider the following mixed $L^1$/$L^2$-optimal control problem.
\begin{problem}[$L^1$/$L^2$-Optimal Control]
\label{prob:L1L2}
Find an admissible control $\{\vc{u}(t): t\in[0,T]\}\in\U$ that minimizes
\begin{equation}
 \begin{split}
  J_{12}(\vc{u}) 
  &\eq 
   \sum_{i=1}^m \biggl(\lambda_i\|u_i\|_{L^1}+\frac{1}{2}r_i\|u_i\|_{L^2}^2\biggr)\\
  &=  
   \int_0^T \sum_{i=1}^m \biggl(\lambda_i |u_i(t)| + \frac{1}{2}r_i|u_i(t)|^2\biggr) dt,
 \end{split}
 \label{eq:J_12}
\end{equation}
where $\lambda_i>0$ and $r_i>0$, $i=1,\dots,m$, are given weights.
\end{problem}

To discuss the optimal solution(s) of the above problem,
we next give necessary conditions for the $L^1$/$L^2$-optimal control
using the minimum principle of Pontryagin.

The Hamiltonian function is given by
\[
 \begin{split}
  H(\vc{x},\vc{p},\vc{u})
   &= \sum_{i=1}^m \biggl(\lambda_i|u_i|+\frac{1}{2}r_i|u_i|^2\biggr)\\
    &\quad + \vc{p}^\top \biggl(\vc{f}(\vc{x})+\sum_{i=1}^m\vc{g}_i(\vc{x})u_i\biggr)
 \end{split}
\]
where $\vc{p}$ is the costate vector.
Let $\vc{u}^\ast$ denote the optimal control and
$\vc{x}^\ast$ and $\vc{p}^\ast$ the resultant optimal state and costate,
respectively.
Then we have the following result.
\begin{lem}
\label{lem:L1L2}
The $i$-th element $u_i^\ast(t)$ of
the $L^1$/$L^2$-optimal control $\vc{u}^\ast(t)$
satisfies
\begin{equation}
 u_i^\ast(t) = -\sat\left\{\shrink_{\lambda_i/r_i}\left(r_i^{-1}\vc{g}_i\bigl(\vc{x}^\ast(t)\bigr)^\top\vc{p}^\ast(t)\right)\right\},
 \label{eq:uopt_1}
\end{equation}
where 
$\shrink_{\lambda/r}(\cdot)$ is the shrinkage function defined by
\[
 \shrink_{\lambda/r}(v)  \eq
   \begin{cases} 
     v+\lambda/r& {\rm{if~}} v<-\lambda/r,\\
     0, & {\rm{if~}} -\lambda/r\leq v \leq \lambda/r,\\
     v-\lambda/r, & {\rm{if~}} \lambda/r<v,\\
   \end{cases}
\]
and
$\sat(\cdot)$ is the saturation function defined by
\[
 \sat(v)  \eq
   \begin{cases} 
     -1, & \text{\rm if~} v<-1,\\
     v, & \text{\rm if~} -1\leq v \leq 1,\\
     1, & \text{\rm if~} 1<v.\\
   \end{cases}
\]
See Figs.~\ref{fig:shrink} and \ref{fig:sat_shrink} for the graphs of
$\shrink_{\lambda/r}(\cdot)$ and 
$\sat\!\left(\shrink_{\lambda/r}(\cdot)\right)$,
respectively.
\begin{figure}[tb]
\centering
\includegraphics[width=0.9\linewidth]{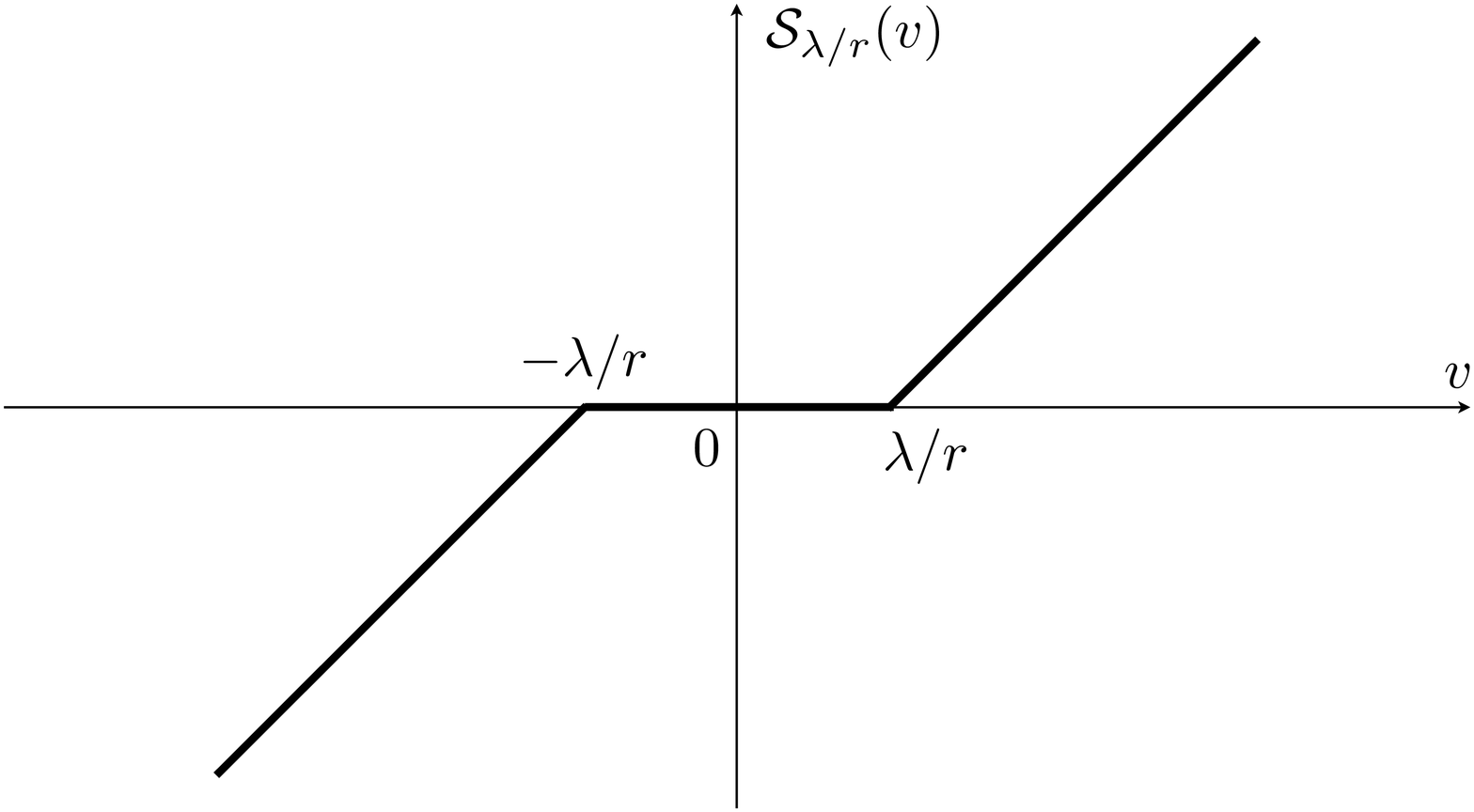}
\caption{Shrinkage function $\shrink_{\lambda/r}(v)$}
\label{fig:shrink}
\end{figure}
\begin{figure}[tb]
\centering
\includegraphics[width=0.9\linewidth]{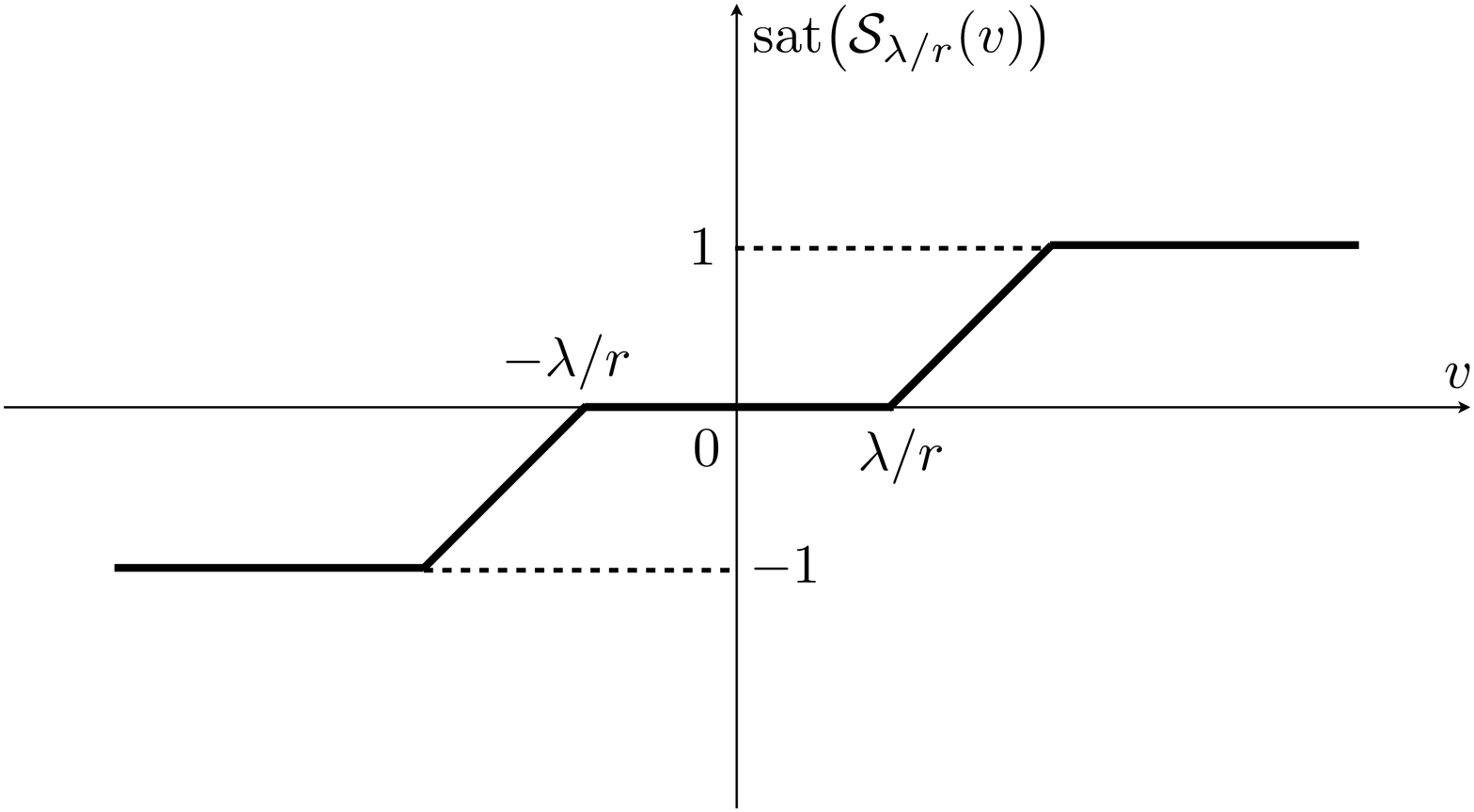}
\caption{Saturated shrinkage function $\sat\!\left(\shrink_{\lambda/r}(v)\right)$}
\label{fig:sat_shrink}
\end{figure}
\end{lem}

{\bf Proof.}
The result is easily obtained from the fact that
\[
 -\sat\left\{\shrink_{\lambda/r}\left(r^{-1}a\right)\right\} = \argmin_{|u|\leq 1} \lambda |u| + \frac{1}{2}r|u|^2+au,
\]
for any $\lambda>0$, $r>0$, and $a\in\R$.
\hfill\QED

From Lemma \ref{lem:L1L2}, we have the following proposition.
\begin{prop}[Continuity]
\label{prop:continuity}
The $L^1$/$L^2$-optimal control $\vc{u}^\ast(t)$ is continuous in $t$ over $[0,T]$.
\end{prop}

{\bf Proof.}
Without loss of generality, we assume $m=1$
(a single input plant),
and omit subscripts for $u$, $r$, $\lambda$, and so on.
Let
\[
 \bar{u}(\vc{x},\vc{p}) \eq -\sat\left\{\shrink_{\lambda/r}\left(r^{-1}\vc{g}_i(\vc{x})^\top\vc{p}\right)\right\}.
\]
Since functions $\left(\sat\circ\shrink_{\lambda/r}\right)(\cdot)$ and $\vc{g}(\cdot)$
are continuous,
$\bar{u}(\vc{x},\vc{p})$ is also continuous in $\vc{x}$ and $\vc{p}$.
It follows from Lemma~\ref{lem:L1L2} that
the optimal control $u^\ast$ given in \eqref{eq:uopt_1}
is continuous in
$\vc{x}^\ast$ and $\vc{p}^\ast$.
Hence, $u^\ast(t)$ is continuous if $\vc{x}^\ast(t)$ and $\vc{p}^\ast(t)$
are continuous in $t$ over $[0,T]$.

The canonical system for the $L^1$/$L^2$-optimal control is given by
\[
 \begin{split}
  \frac{d\vc{x}^\ast(t)}{dt} &=  \vc{f}\bigl(\vc{x}^\ast(t)\bigr) + \vc{g}\bigl(\vc{x}^\ast(t)\bigr)\bar{u}\bigl(\vc{x}^\ast(t),\vc{p}^\ast(t)\bigr),\\
  \frac{d\vc{p}^\ast(t)}{dt} &= - \vc{f}'\bigl(\vc{x}^\ast(t)\bigr)^\top \vc{p}^\ast(t)\\
  	&\qquad - \bar{u}\bigl(\vc{x}^\ast(t),\vc{p}^\ast(t)\bigr) \vc{g}'\bigl(\vc{x}^\ast(t)\bigr)^\top \vc{p}^\ast(t).
 \end{split} 
\]
Since $\vc{f}(\vc{x})$, $\vc{g}(\vc{x})$, $\vc{f}'(\vc{x})$, and $\vc{g}'(\vc{x})$ are continuous
in $\vc{x}$ by assumption,
and so is $\bar{u}(\vc{x},\vc{p})$ in $\vc{x}$ and $\vc{p}$,
the right hand side of the canonical system is continuous in $\vc{x}^\ast$ and $\vc{p}^\ast$.
From a continuity theorem of dynamical systems,
e.g. \cite[Theorem 3-14]{AthFal},
it follows that the resultant trajectories $\vc{x}^\ast(t)$ and $\vc{p}^\ast(t)$ are
continuous in $t$ over $[0,T]$.
\hfill\QED

Proposition \ref{prop:continuity} motivates us to use $L^1/L^2$ optimization
as in Problem \ref{prob:L1L2} for continuous hands-off control.

In general, the degree of continuity (or smoothness) and the sparsity of the control input
cannot be optimized at the same time.
Then, the weights $\lambda_i$ or $r_i$ can be used for
trading smoothness for sparsity.
From Lemma~\ref{lem:L1L2},
increasing the weight $\lambda_i$ (or decreasing $r_i$)
makes the $i$-th input $u_i(t)$ sparser
(see also Fig.~\ref{fig:sat_shrink}).
On the other hand,
decreasing $\lambda_i$ (or increasing $r_i$)
smoothens $u_i(t)$.
Moreover, we have the following limiting properties.
\begin{prop}[Limiting property]
Assume the $L^1$-optimal control problem 
is normal.
Let $\vc{u}_1(\vc{\lambda})$ and $\vc{u}_{12}(\vc{\lambda},\vc{r})$
be solutions to respectively Problems \ref{prob:L1} and \ref{prob:L1L2}
with parameters
\[
 \vc{\lambda}\eq(\lambda_1,\dots,\lambda_m),\quad
 \vc{r}\eq(r_1,\dots,r_m).
\]
For any fixed $\vc{\lambda}>0$,
we have
\[
 \lim_{\vc{r}\rightarrow\vc{0}} \vc{u}_{12}(\vc{\lambda},\vc{r}) = \vc{u}_1(\vc{\lambda}).
\] 
Moreover, for any fixed $\vc{r}>0$,
we have
\[
 \lim_{\vc{\lambda}\rightarrow\vc{0}}\vc{u}_{12}(\vc{\lambda},\vc{r}) = \vc{u}_2(\vc{r}),
\]
where $\vc{u}_2(\vc{r})$ is an $L^2$-optimal (or minimum-energy) control
discussed in \cite[Chap.~6]{AthFal},
that is,
a solution to a control problem where $J_1(\vc{u})$ in Problem \ref{prob:L1} is replaced with
\begin{equation}
 J_2(\vc{u}) = \int_0^T \sum_{i=1}^m \frac{1}{2}r_i|u_i(t)|^2dt.
 \label{eq:J2}
\end{equation}
\end{prop}

{\bf Proof.}
The first statement follows directly from the fact that
for any fixed $\lambda>0$, we have
\[
 \lim_{r\rightarrow 0}\sat\!\left({\mathcal{S}}_{\lambda/r}(r^{-1}w)\right) = \dez_{\lambda}(w),\quad
 \forall w\in\R\setminus\{\pm \lambda\},
\]
where $\dez_\lambda(\cdot)$ is the dead-zone function defined in \eqref{eq:dez}. 
The second statement derives from the fact that
for any fixed $r>0$, we have
\[
 \lim_{\lambda\rightarrow 0}\sat\!\left({\mathcal{S}}_{\lambda/r}(v)\right) = \sat(v),\quad
 \forall v\in\R.
\]
\hfill\QED

In summary, the $L^1$/$L^2$-optimal control is an \emph{intermediate control} between 
the $L^1$-optimal control (or the maximum-hands-off control)
and the $L^2$-optimal control.

\section{Examples}
\label{sec:examples}
We here consider the following 4th order system:
\[
 \frac{d\vc{x}(t)}{dt} =
  \begin{bmatrix}0&-1&0&0\\1&0&0&0\\0&1&0&0\\0&0&1&0\end{bmatrix}\vc{x}(t)
  + \begin{bmatrix}2\\0\\0\\0\end{bmatrix}u(t).
\]
We set the final time $T=10$, and the initial and final states as
\[
 \vc{x}(0) = [1,1,1,1]^\top,\quad \vc{x}(10)=\vc{0}.
\] 
Note that the system has poles at $s=0$, $0$, $\pm j$.

\begin{figure}[tb]
\centering
\includegraphics[width=\linewidth]{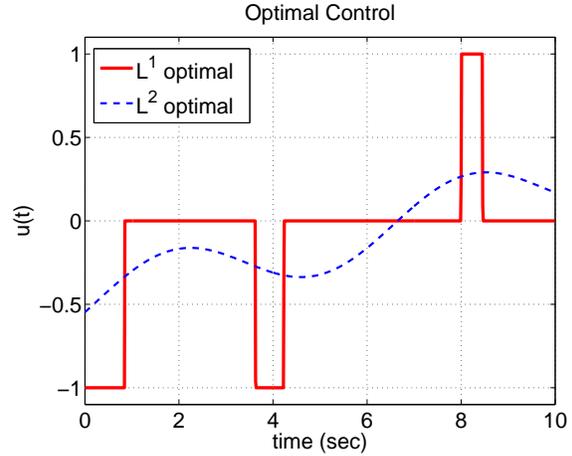}
\caption{Maximum-hands-off control via $L^1$ optimization (solid) and $L^2$ optimal control (dashed).}
\label{fig:control_4th}
\end{figure}

\begin{figure}[tb]
\centering
\includegraphics[width=\linewidth]{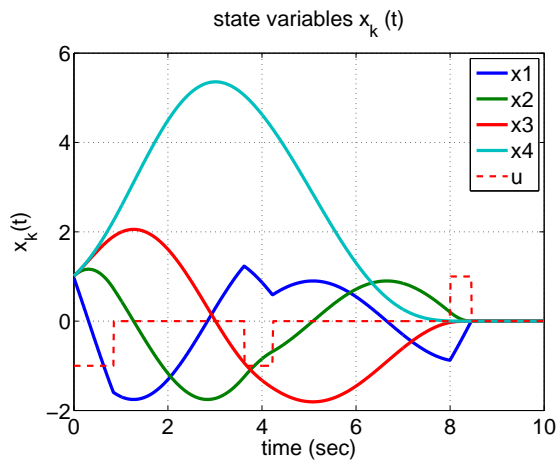}
\caption{Maximum-hands-off control:
state variables $x_1(t),\dots,x_4(t)$ (solid) and input $u(t)$ (dashed)}
\label{fig:states_4th}
\end{figure}

We first compute the maximum-hands-off control
with $L^1$-optimal control as discussed in Section \ref{sec:main}.
We compute the optimal control input by a time discretization method,
see e.g., \cite[Sec.~2.3]{Ste}.
Fig.~\ref{fig:control_4th} shows the obtained control.
The figure also shows the $L^2$-optimal control that minimizes $J_2(u)$
in \eqref{eq:J2} with $r_1=1$.
We can see that the maximum-hands-off control is quite sparse.
In fact, we have
\[
\|u\|_{L^0} = 1.92 \text{~(sec)},
\]
which is $19.2$\% out of $10$ (sec).
In other words, the control keeps hands-off over $80.8$\% of the control period.
On the other hand, the $L^2$ optimal control is not sparse,
while its energy, $J_2(u)$,  is smaller than that of maximum-hands-off control.
Fig.~\ref{fig:states_4th} shows the state variables $x_1(t)$, $x_2(t)$,
$x_3(t)$, and $x_4(t)$
along with the maximum-hands-off control $u(t)$ over time interval $[0,10]$.
We can see that the states almost stay at the origin
after the last switching time, $t=8.47$ (sec).

\begin{figure}[t!]
\centering
\includegraphics[width=\linewidth]{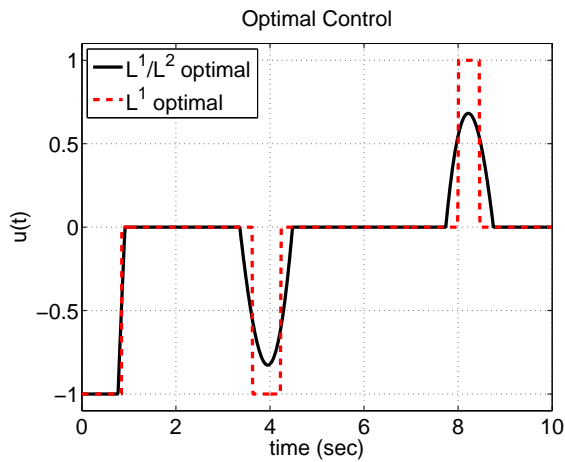}
\caption{$L^1$/$L^2$-optimal control (solid) and maximum-hands-off control (dash).}
\label{fig:control_L1L2}
\end{figure}

\begin{figure}[t!]
\centering
\includegraphics[width=\linewidth]{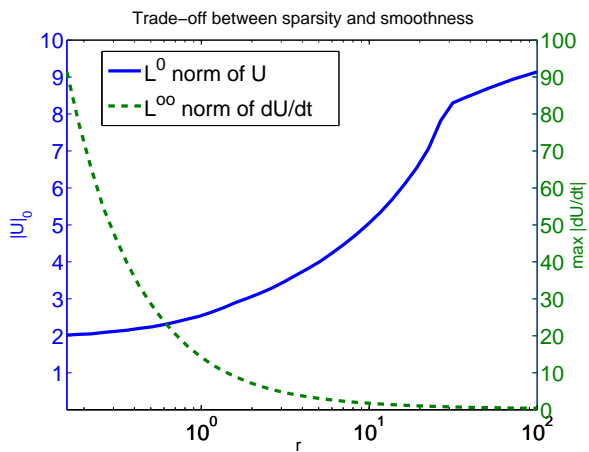}
\caption{$L^0$ norm of $u(t)$ (solid) and $L^\infty$ norm of $du/dt$ (dash) versus weight $r$.}
\label{fig:tradeoff}
\end{figure}

We next consider the $L^1$/$L^2$-optimal control method
proposed in Section~\ref{sec:L1L2}.
We use the same parameters as above.
The weights $\lambda_1$ and $r_1$ in \eqref{eq:J_12} are chosen as $\lambda_1=r_1=1$.
We solve the optimal control problem via a time discretization method.
Fig.~\ref{fig:control_L1L2} shows the obtained $L^1$/$L^2$-optimal control.
The figure also shows the maximum-hands-off control
obtained above.
We can see that the $L^1$/$L^2$-optimal control is continuous while
the maximum-hands-off control exhibits
the "bang-off-bang" property.
On the other hand, the $L^1$/$L^2$-optimal control
has a longer support than the $L^1$-optimal control.
To see the tradeoff property between sparsity and smoothness of control,
we compute the $L^0$ norm, $\|u\|_{L^0}$,
and the $L^\infty$ norm of the derivative of $u(t)$, that is,
\[
 \left\|\frac{du(t)}{dt}\right\|_{L^\infty} \eq \sup_{t\in[0,T]} \left|\frac{du(t)}{dt}\right|,
\]
as a function of $r_1$ while $\lambda_1$ is fixed to be $1$.
Fig.~\ref{fig:tradeoff} shows the result.
We can see that the weight $r$ can take account of the tradeoff
between sparsity and smoothness in hands-off control.

\section{Conclusion}
\label{sec:conclusion}
In this article, we have 
presented maximum-hands-off control and shown that it is $L^1$ optimal.
This shows that efficient optimization methods for $L^1$ problems can be used to obtain maximum-hands-off control.
We have also proposed an $L^1$/$L^2$-optimal control to obtain
smooth hands-off control, while the maximum-hands-off control
is discontinuous due to the "bang-off-bang" property.
Numerical examples show the effectiveness of the proposed control.
Future work may include adaptation of hands-off control to sparsely packetized predictive control
as in \cite{NagQue11,NagQueOst12b}.



\end{document}